\documentclass{amsart}
\usepackage{amssymb,amsmath,amsthm}
\usepackage{enumerate}
\usepackage{eucal}
\usepackage{mathabx}
\usepackage{xcolor}
\usepackage{tasks}
\usepackage{graphicx}
\usepackage{hyperref}
\usepackage{mathtools}
\usepackage{etoolbox}
\usepackage{orcidlink}

\hypersetup{
    colorlinks,
    linkcolor={red!75!black},
    citecolor={green!65!black},
    urlcolor={magenta!75!black}
}

\usepackage[backend=biber,sorting=nyt,style=alphabetic,minalphanames=3,maxbibnames=6]{biblatex}

\DeclareNameAlias{author}{last-first}
\theoremstyle{definition}

\theoremstyle{definition}
\newtheorem{Problem}{Problem}

\usepackage{graphicx}  

\title{The number of countable models of first-order theories}
\author[A. Pillay]{Anand Pillay \orcidlink{0000-0002-2471-4097}}

\address[A. Pillay]{University of Notre Dame, Notre Dame, Indiana}

\email[A. Pillay]{apillay@nd.edu}

\thanks{The first author was supported by  NSF grants DMS-2054271 and DMS-2502292 }

\author[P.\ Tanovi\'c]{Predrag Tanovi\'c \orcidlink{0000-0003-0307-7508}}

\address[P.\ Tanovi\'c]{Mathematical Institute SANU\\ Knez Mihailova 36, Belgrade, Serbia}

\email[P.\ Tanovi\'c]{tane@mi.sanu.ac.rs}

\thanks{The second author was supported by the Ministry of Education, Science and Technological Development of Serbia through Mathematical Institute SANU.}

\begin{document}

\begin{abstract} Let $T$ be a complete first-order theory in a countable language that has infinite models. 
We discuss and survey the work done on the number of nonisomorphic countable models $T$ from the perspective of first-order model theory: 
Vaught's conjecture, Martin's conjecture, and problems related to Ehrenfeuch theories (theories with more than one but only finitely many countable models). 
\end{abstract}

\maketitle

Throughout, $T$ denotes a complete first-order theory in a countable language $L$ that has infinite models and $I(\aleph_0,T)$ denotes the number of countable models of $T$, up to an isomorphism. 
To determine $I(\aleph_0,T)$, it suffices to consider only countable models of $T$ with domain $\omega$; since there are at most continuum many $L$-structures with domain $\omega$, $I(\aleph_0,T)\leqslant 2^{\aleph_0}$ holds. Theories with $I(\aleph_0,T)=1$ are the $\aleph_0$-categorical theories. These include the theory of an infinite set, the theories of infinite-dimensional vector spaces over a finite field, and the theory of dense linear orders. Theories with $I(\aleph_0,T)<2^{\aleph_0}$ are said to have few countable models. 

\section{Vaught's conjecture}

The number of countable models was first studied by Vaught in his 1959 conference talk and the subsequent paper \cite{Vaught1961}. He proved that $I(\aleph_0,T)$ is never equal to two and, at the end of his talk, posed a problem: {\it Can it be proved, without the use of the continuum hypothesis, that there exists a complete theory having exactly $\aleph_1$ non-isomorphic denumerable models?} What is commonly known as Vaught's conjecture is a statement predicting a negative answer to Vaught's question:

\smallskip
{\bf Vaught's conjecture.} \ $I(\aleph_0,T)$  is countable or equal to $2^{\aleph_0}$.

\smallskip
Vaught's conjecture naturally extends to or generalizes to classes of countable models defined by a sentence of the infinitary logic $L_{\omega_1,\omega}$, which permits countable conjunctions and
disjunctions, but allows only finite strings of quantifiers: every $L_{\omega_1,\omega}$-sentence has either countably many or $2^{\aleph_0}$ countable models. 
Morley in \cite{Morley1970number} made the first progress toward Vaught's conjecture: he proved that the only infinite values possible for the number of countable models of an $L_{\omega_1,\omega}$-sentence are $\aleph_0,\aleph_1$, and $2^{\aleph_0}$. This is still the strongest result known for an arbitrary first-order theory; it will be discussed in more detail in Section \ref{Section_martin}. 
Morley's theorem turns out to be rather a result in descriptive set theory; it is a special case of Burgess's theorem on analytic equivalence relations \cite{Burgess1978equivalences}. Statements equivalent to Vaught's conjecture were found in descriptive set theory (Steel \cite{Steel1978vaught}), topological dynamics (Becker and Kechris \cite{Becker-Kechris}) and computability theory (Sacks \cite{SACKS1983185} and Montalban \cite{Montalban2013computability}). 
Notice that Vaught's conjecture, as stated, is trivially true under the Continuum Hypothesis (CH). However, there are statements provably equivalent in ZF+$\lnot$CH to Vaught's conjecture that do not follow from CH. Steel found one in \cite{Steel1978vaught}, which he called the strong Vaught conjecture.

\begin{itemize}\item[{\bf VC1}]\  Every $L_{\omega_1,\omega}$-sentence $\phi$ has either countably many models or there is a perfect subset $P$ of $^{\omega}2$  such that any two distinct elements of $P$ code nonisomorphic models of $\phi$.
\end{itemize}

Another statement equivalent to Vaught's conjecture derives from Morley's work \cite{Morley1970number}, in which a modification of Scott's analysis of the isomorphism relation was used; Morley's work will be further discussed in Section \ref{Section_martin}. Scott in \cite{Scott2014logic} proved that for every countable model $M$ there exists an $L_{\omega_1,\omega}$-sentence $\phi_M$ such that the countable models of $\phi_M$ are exactly those that are isomorphic to $M$; this sentence is called a Scott sentence for $M$. The smallest possible ordinal value of the quantifier rank of a Scott sentence for $M$ is called the Scott rank of $M$. The Scott rank of the theory, or an $L_{\omega_1,\omega}$-sentence, is the supremum of the Scott ranks of its countable models.   

\begin{itemize}
    \item[{\bf VC2}] Every $L_{\omega_1,\omega}$-sentence with few countable models has countable Scott rank. 
\end{itemize}
  
Vaught's conjecture is widely open in both directions. By now, it has been confirmed for various special classes of theories, mostly using first-order model theory techniques. The only exception is the class of theories of colored trees, for which Steel in \cite{Steel1978vaught} proved {\bf VC1} by using descriptive set theory and infinitary model theory methods; so far, no first-order proof has been published.   

So, here are some classes of first-order theories for which Vaught's conjecture has been proved. Some other very special classes will be mentioned later. 
\begin{itemize}
    \item Theories of colored orders (Rubin in \cite{rubin1974theories};
    \item Theories of one unary operation (Miller in \cite{miller1981vaught});
    \item $\omega$-stable theories (Shelah, Harrington and Makkai in \cite{Shelah1984proof};
    \item Stable theories with Skolem functions (Lascar in \cite{Lascar1981skolem});
     \item o-minimal theories (Mayer in \cite{Mayer1988vaught});    
    \item Varieties (equational theories) (Hart, Starchenko and Valeriote in \cite{HartStarchenkoVal1994vaught});
    \item Superstable theories of finite U-rank (Buechler in \cite{Buechler2008vaught}).
\end{itemize}
In all the cases listed, a purely first-order model-theoretic approach to Vaught's conjecture via the following broad scheme is used: restrict to a class of theories, isolate the properties of theories that cause $I(\aleph_0,T)=2^{\aleph_0}$, and then,  
 assuming that $T$ does not have any of them, find a reasonable system of invariants for isomorphism types of countable models of $T$. 
Here, from our point of view, theories with $I(\aleph_0,T)=2^{\aleph_0}$ are considered to have the isomorphism problem maximally complicated, as any further distinction between them inevitably leads to invariant descriptive set theory; see  \cite{Friedman1989borel}, \cite{Marker2007borel}, \cite{laskowski2015borel}, and \cite{laskowski2023most}. 
By reasonable invariants, structures that can be easily counted are usually meant, such as cardinals or trees of cardinals with fixed finite height. However, invariants described by $L_{\omega_1,\omega}$-sentences of countably bounded quantifier rank are good enough to confirm {\bf VC2}. For example, the invariants obtained for the $\omega$-stable theories in \cite{Shelah1984proof} were originally described by $L_{\omega_1,\omega}$-sentences of quantifier rank $<\omega\cdot\omega$; later, Bouscaren in \cite{Bouscaren1984martin} reduced the bound to $\omega+\omega$. 
 
One of the main challenges in classification is establishing enough conditions that lead to $I(\aleph_0,T)=2^{\aleph_0}$, thus enabling the identification of relevant invariants. Vaught noted one such general condition: $T$ is not {\em small} (there are uncountable many pure types without parameters).
In all the cases listed, these conditions are quite different and no new general condition implying
$I(\aleph_0,T)=2^{\aleph_0}$ was found. Therefore, it is unlikely to expect that the confirmation of Vaught's conjecture can come from first-order model theory. It is more likely that studying more complicated classes of theories may lead to a possible counterexample. 
Nevertheless, the procedure of ``working from the bottom up", is interesting in its own right, as more information is given about describing the countable models of specific kinds of theories. 

Conditions that lead to $I(\aleph_0,T)=2^{\aleph_0}$ are usually model-theoretical in nature and, in some cases, verifying these conditions can be difficult, even within certain specific classes of algebraic structures. An example of this is the theory of differentially closed fields of characteristic 0. The theory DCF$_0$ has been known to be $\omega$-stable and, from Shelah's proof of Vaught’s conjecture for $\omega$-stable theories, the number of its countable models is $\aleph_0$ or $2^{\aleph_0}$. However, it took over a decade before Hrushovski and Sokolovi\'c (unpublished) proved that DCF$_0$ satisfies one of the conditions leading to $2^{\aleph_0}$ models (ENI-DOP). Their proof made use of the truth of the Zilber dichotomy for strongly minimal sets in DCF$_0$. Another proof appears in \cite{Pillay2005model}.

\smallskip 
Let us discuss in more detail both existing results and conjectures in some broad classes of theories.

\medskip
 (1) {\bf Stable theories.} 
 In terms of studying the number of {\em uncountable} models of a countable complete first-order theory $T$, Shelah's program \cite{Shelah1990classification} was very successful. His Main Gap Theorem states that either $I(\kappa,T) = 2^{\kappa}$ for all uncountable $\kappa$ or $I(\aleph_{\alpha}, T)$ is less than $\beth_{\omega_{1}}(|\omega + \alpha|)$ for all $\alpha\geq 1$, where of course $I(\kappa, T)$ is the number of models of $T$ of cardinality $\kappa$.  The notion of stability is crucial: $T$ (countable complete) is stable if for some infinite $\kappa$ the number of complete $1$-types over any model of cardinality $\kappa$ is at most $\kappa$.  If $T$ is unstable then it fits into the first possibility above, the maximum number of models in any uncountable cardinality. So, Shelah could restrict attention to stable theories, for which either one again gets the maximum number of models or obtains some kind of classification or description of models (leading to the second possibility above). An enormous machinery was developed for understanding and describing models of stable theories as well as definability in such models. However, on the face of it, the methods did not apply to the case of countable models.

 Shelah built on the fundamental and early work of Morley, stating that $I(\kappa,T) = 1$ for some uncountable $\kappa$ iff $I(\kappa,T) = 1$ for all uncountable $T$ (these are the so-called uncountably categorical theories).  Morley's work also showed that for such theories $I(\omega, T) \leq \omega$ (establishing Vaught's conjecture for uncountably categorical $T$).  But in \cite{Baldwin-Lachlan} it was shown that to {\em any} model $M$ of an uncountably categorical $T$ one can attach a cardinal $dim(M)$ which is determined by and determines $M$ up to isomorphism  and such that $|M| = dim(M) + \omega$, which led to the theorem that such $T$ are either $\omega$-categorical or have exactly $\omega$-many countable models (this is connected with Section 3). A key role (in the Baldwin-Lachlan approach) is given by {\em strongly minimal formulas or definable sets $X$} meaning $X$ is infinite and cannot be partitioned into two infinite definable sets. The unique nonalgebraic type containing a strongly minimal formula is called a strongly minimal type.  The number of (algebraically) independent realizations of $p$ in a model $M$, is well-defined and coincides with $dim(M)$ above. (There are some subtle issues regarding parameters over which the strongly minimal set $X$ is defined.)

 In \cite{Shelah1984proof} the authors were able to adapt Shelah's machinery for studying uncountable models to the countable case, for $\omega$-stable theories. These are countable complete theories which are $\kappa$-stable for all infinite $\kappa$  (as is the case for uncountably categorical $T$). Now $\omega$-stability has many implications, such as the existence of prime models over all sets, extraction of indiscernibles, and ``finite multiplicity" of all types. 
 (The multiplicity of a type $p\in S_n(A)$ is the number of strong types over $A$ that extend $p$; it is either finite or $2^{\aleph_0}$. Two tuples $\bar a,\bar b$ of length $n$ have the same strong type over $A$ if and only if $E(\bar a,\bar b)$ holds for $A$-definable finite equivalence relations on $M^n$; strong types over a finite set $A$ are precisely complete types over $acl^{eq}(A)$.) 
 The ubiquity of {\em strongly regular} (SR) types which support well-behaved dimensions is crucial. (Strongly regular types generalize strongly minimal types.) A complete type over a finite set is called eventually non-isolated (eni) if some nonforking extension over a larger finite set is non-isolated; an ENI type is a strongly regular eni type. Under the assumption that $T$ has fewer than continuum many countable models, it is proved that $T$ has the technical property ``ENI depth at most $1$". From this a structure theorem for countable models of $T$ is proved (involving dimensions of suitable strongly regular types) from which it is concluded that $T$ has either countably many or continuum many countable models.

 The next step in the stability hierarchy was the case of superstable theories, which is still open in full generality. 
 As for the purposes of looking at Vaught's conjecture, one may assume that $T$ is {\em small}  (countably many complete types without parameters), the only difference with the $\omega$-stable case is that complete types over finite sets may have infinite multiplicity. Attention was first focused on the case where $T$ has ``finite rank" and even when $T$ is {\em weakly minimal}, meaning that there are a bounded number (at most continuum) of nonalgebraic $1$-types (over any model). In the work of Saffe (unpublished) on the problem, a certain proposition was claimed; that a rank $1$ (i.e. weakly minimal) complete type $p$ over a finite set $A$, which has infinite multiplicity, must be isolated, unless $T$ has continuum many countable models.  A mistake was found in the proof. Subsequently, Buechler \cite{Buechler1985classification} gave a proof of Vaught's conjecture for weakly minimal theories, assuming the truth of ``Saffe's hypothesis". Another ingredient came from ``geometric stability theory": a rank $1$ type is either locally modular or has Morley rank $1$. (This is Buechler's dichotomy). 
 Then in a breakthrough result, Newelski \cite{Newelski1990saffe} gave a proof of Saffe's hypothesis. Newelski's proof introduced new techniques for dealing with and understanding strong types, and led to a study of ``small profinite structures and groups".  An account of Newelski's work and the proof of Vaught's conjecture for weakly minimal theories (under a convenient additional assumption of ``unidimensionality") appears in Chapter 6 of \cite{Pillay1996geometric}. 

 The main question at that point (around 1990) was to try to extend the proof of the weakly minimal case to general superstable theories.  The superstable finite rank case was eventually proved by Buechler \cite{Buechler2008vaught}, appearing considerably later and will be discussed more below. But in the meantime, there were some other observations made, also with a geometric-stability-theoretic flavor. In \cite{Low-Pillay1992superstable}, Low and the first author proved that if $T$ is (countable, complete)  superstable, not $\omega$-stable and with fewer than continuum many countable models, then there is a definable group $G$ with locally modular regular generics such that any nongeneric types in $G$ are $\omega$-stable (i.e. have Morley rank).  One conclusion was that a superstable theory with no infinite definable groups satisfies Vaught's conjecture.  It was also natural to study objects given by the conclusion of the theorem above, as structures in their own right.  Namely, to consider the theory of a locally modular regular non connected-by-finite abelian group $G$ such that the nongeneric part of $G$ is $\omega$-stable.  The first author did precisely this in \cite{Pillay1992certain} (unpublished) proving Vaught's conjecture by adapting ideas from the weakly minimal case.  
 For the general superstable (possibly infinite rank) case one would also like some version of the Buechler dichotomy for infinite rank regular types, such as strongly regular or locally modular. 

 We return to Buechler's solution to Vaught's conjecture in the superstable, finite $U$-rank setting \cite{Buechler2008vaught}.  This was a tour-de-force. The so-called semiminimal analyses of elements and types were crucial, but also some representation theory played a role. One would have expected the paper to be picked up by the community as the first step in an inductive proof of Vaught's conjecture for superstable theories.  But for various reasons (maybe the extremely technical nature of the proof,  the delay in coming to publication, or even the  many additional levels of complexity compared with the weakly minimal case), the actual and unfortunate consequence was the general cessation of work on the subject.  It would be desirable to have a clear conceptual account of this superstable finite rank case of Vaught's conjecture, maybe under an additional assumption of ``unidimensionality". Unidimensionality means that any two types have some kind of interaction (over additional parameters) and focusing on this case concentrates on the truly new phenomena one must deal with. 
 
 Nevertheless, some earlier work dealt with very special cases, such as making the assumption of $1$-based. This is a kind of global local modularity, specifically that the canonical base of a stationary type $tp(a/A)$ is in the algebraic closure of $a$ (in $T^{eq}$).  There are many consequences of $1$-basedness including the structure of definable groups. 
 The second author, in \cite{Tanovic2007nonisol}  proved that every superstable one-based theory with few models has $U$-rank strictly below $\omega^\omega$. This is relevant to the understanding of higher-rank regular types. As shown in \cite{Tanovic2013simple}, the same bound would apply without the assumed one-basedness, provided that the answer to the following problem from geometric model theory is affirmative. 

 \begin{Problem}
The generic type of any simple superstable group is eventually non-isolated.
 \end{Problem}

 Among $1$-based theories are complete theories of $R$-modules for $R$ a (countable) ring. All such theories are stable, and there is a nice characterization of the superstable such theories. Rather surprisingly (bearing in mind the fairly comprehensive understanding of the class of models), the following is still open:

 \begin{Problem} 
 Prove Vaught's conjecture for superstable theories of modules.   
 \end{Problem}

 Let us finally consider the case of general stable theories. As in the case of Problem 2 above, Vaught's conjecture remains open for arbitrary complete theories of modules (over a countable ring). Additional assumptions may give ``superstable-like" behaviour. 
  For example in \cite{Herwig1992stable} it is shown that under an additional assumption of smallness (countably many complete types over $\emptyset$) types have finite weight and there are enough regular types. 
  Under a $1$-based assumption, the technical notion of ``triviality" is equivalent to the nonexistence of definable infinite groups. The following was proved in \cite{Pillay1992countable}  (making use also of \cite{Low-Pillay1992superstable});  if $T$ is a trivial $1$-based stable theory which is not $\omega$-stable then $T$ has continuum many countable models.

\medskip
(2) {\bf Weakly o-minimal theories.} 
Recall that an $o$-minimal theory is one such that every model is totally ordered by a fixed relation symbol in the language and every definable subset of $1$-space is a finite union of points and intervals (with endpoints in the model). As mentioned earlier, Mayer proved Vaught's conjecture (actually much more) for countable (complete) $o$-minimal theories. 
Weakly minimal theories are a generalization: every definable subset of $1$-space is a finite union of convex sets. Countable models of weakly o-minimal theories are substantially more complicated than those of the o-minimal ones, even in the $\aleph_0$-categorical case. Pillay and Steinhorn classified all $\aleph_0$-categorical o-minimal theories in \cite{PillaySteinhorn1986definableI}; in particular, they are all binary (every formula is equivalent modulo $T$ to a Boolean combination of formulae with at most two free variables). Binary $\aleph_0$-categorical weakly o-minimal theories were classified by Kulpeshov in \cite{kulpeshov2007criterion} and \cite{Kulpeshov2016countably}. However, Herwig et al. in \cite{Herwig-Macpherson1999} suggested that there is no close description in the general weakly o-minimal case; they describe an iterative process for the construction of more and more complicated such structures. This also shows that the proof of Vaught's conjecture in this case would be substantially more complicated than in the o-minimal case.
Some progress has been made recently in \cite{Moconja2020stationarily}, where Vaught's conjecture was confirmed for a subclass of binary linearly ordered theories which includes all weakly o-minimal and all weakly quasi-o-minimal theories (every definable subset of a saturated $M\models T$ is a Boolean combination of unary $L$-definable sets and convex sets), and by Kulpeshov in  \cite{Kulpeshov2020}, who proved Vaught's conjecture for weakly o-minimal theories of finite convexity rank by reducing it to the binary case. 
Recently, in \cite{Moconja2025wom}, the notion of weak o-minimality was introduced for complete types in an arbitrary first-order theory and it is shown that a theory is weakly quasi-o-minimal if and only if all types in $S_1(T)$ are weakly o-minimal. It is shown that these types satisfy a weak version of the Monotonicity Theorem and share some nice properties of weight-one types in stable theories. This suggests that the weakly quasi-o-minimal case of the conjecture is not substantially more complicated than the weakly o-minimal one. 
 
\medskip 
(3) {\bf Theories with built-in Skolem functions.} 
Recall that $T$ is said to have (built-in) Skolem functions if for every formula $\phi(x,y)$ of the language (where $x$, $y$ may be tuples of variables) there is a partial $\emptyset$-definable function $f(y)$ such that in any model of $T$ the following holds: 
$\forall y (\exists x \phi(x,y) \to \phi(f(y), y)))$. A consequence (in fact, an equivalent statement) is that any definably closed subset of a model of $T$ is an elementary substructure. 
The only known confirmed subcases of Vaught's conjecture for such theories are for theories that are also stable (Lascar \cite{Lascar1981skolem}) or admit a definable linear order (Shelah in \cite{Shelah1978endext}). In the latter case, Shelah proved $I(\aleph_0,T)=2^{\aleph_0}$; this can be used in the proof that $I(\aleph_0,T)=2^{\aleph_0}$ holds for every complete first-order $T$ that interprets an infinite discrete linear order (\cite{Tanovic2024vaught}). This suggests that the answer to the following question should be positive, at least for theories with the strict order property; the existence of a definable partial ordering with some infinite chains. 
    
\begin{Problem}
Does every unstable theory with built-in Skolem functions have $2^{\aleph_0}$ countable models?   
\end{Problem}

\section{Martin's conjecture}\label{Section_martin}

Martin's conjecture is a strengthening of Vaught's conjecture. It predicts certain $L_{\omega_1,\omega}$-sentences of quantifier rank $\leqslant \omega+\omega$ as invariants of the isomorphism types of countable models of theories with few countable models. The conjecture was motivated by the proof of Morley's theorem from \cite{Morley1970number}, where a modification of Scott's analysis of isomorphism types of countable models was used. 
First, Morley proved that every theory $T$ (or even a $L_{\omega_1,\omega}$-sentence) with few countable models is scattered, that is, for every countable $L_{\omega_1,\omega}$-fragment, $\mathcal F$, there are only countably many $\mathcal F$-types consistent with $T$. Next, assuming that $T$ is scattered Morley defines $L_1(T)$, the smallest (countable) fragment of $L_{\omega_1,\omega}$ that contains all first-order formulas and formulae $\bigwedge_{\phi(x)\in p}\phi(x)$ for all $p\in S_n(T)$ ($n\in\mathbb N$),  
and a continuous sequence of countable fragments, $(L_\alpha(T)\mid 1\leqslant \alpha<\omega_1)$, such that $L_{\alpha+1}(T)$ is the smallest fragment that contains $L_{\alpha}(T)$ and all formulae that express the conjunction of all formulae of some complete $L_\alpha(T)$-type. Then he proved that any countable model has a Scott sentence within $\bigcup_{\alpha<\omega_1}L_{\alpha}(T)$. Since each $L_\alpha(T)$ is countable, this leads to the conclusion that either there exists some $\alpha<\omega_1$ such that every countable model has a Scott sentence in $L_\alpha(T)$ - thus yielding $I(T,\aleph_0 )\leqslant \aleph_0$ - or alternatively, $I(T,\aleph_0)=\aleph_1$. Martin's conjecture predicts the complete $L_1(T)$-theories as
invariants of the isomorphism types of countable models of $T$ so that, in particular, Scott sentences can be found even in $L_2(T)$. For each $M\models T$ let $T_1(M)$ denote the $L_1(T)$-theory of $M$; notice that the quantifier rank of $T_1(M)$ is at most $\omega+\omega$.  

\smallskip
 {\bf Martin's Conjecture.} If $T$ has few countable models, then $T_1(M)$ is $\aleph_0$-categorical for all countable $M\models T$.

\smallskip
  
Among the $\mathcal L_{\omega_1,\omega}$-fragments, the fragment $L_1(T)$ is quite reasonable from the first-order perspective, since (roughly speaking) its sentences describe finitary relations between the loci of types in a model. This is why Martin's conjecture is a matter of first-order rather than infinitary model theory. However, it is also reasonable to try finding a counterexample to Martin's conjecture as it may be a starting point in constructing a counterexample to Vaught's conjecture. 

For each aforementioned class of theories, the classification problem has been explicitly solved, or its solution can be derived from the proof. The invariants obtained for each class are quite simple, so the validity of Martin's conjecture was derived either as a straightforward corollary (in \cite{Mayer1988vaught} and \cite{Buechler2008vaught}), or by using some mild additional arguments: C.M.Wagner in \cite{CMWagner1982martin} for theories of colored orders and theories of one unary operation, Bouscaren in \cite{Bouscaren1984martin} for $\omega$-stable theories. 
The only class of theories for which the Vaught conjecture has been resolved but Martin'sconjecture is still open is that of theories of colored trees. Steel in \cite{Steel1978vaught} used descriptive set theory to prove the version {\bf VC1} for colored trees, but it is unlikely that his arguments can be adapted to yield invariants sufficiently simple to confirm Martin's conjecture; C.M.Wagner in \cite{CMWagner1982martin} derives the bound $\omega^2$ for the Scott rank.

\begin{Problem}
Prove Martin's conjecture for theories of colored trees.    
\end{Problem}

\section{Ehrenfeucht theories}

 Theories with $1<I(\aleph_0,T)<\aleph_0$ are called Ehrenfeucht theories. These were named after Ehrenfeucht, who gave the first example of a theory with $I(\aleph_0,T)=3$, which is the theory of a dense linear order without endpoints equipped with an increasing sequence of constants: the language has a binary relation symbol $<$ and infinitely many constant symbols $c_0,c_1,...$ with axioms stating that $<$ is a dense linear order without endpoints with $c_0<c_1<c_2<...$. The three models are determined by whether this sequence is unbounded, has a supremum, or is bounded but with no supremum. It is easy to modify this example to obtain a theory with $I(\aleph_0,T)=n$ for any $n>3$ in the following way: Add to the language $n-2$ unary predicates (colours) $P_1,...P_{n-2}$ with axioms that state that the colours partition the domain, that elements of each colour form a dense subset, and that $P_1(c_i)$ holds for all $i\in\omega$. The countable models for which $\lim c_i=c$ are divided into $n-2$ isomorphism types depending on the color of $c$. 

Another example of a theory with $I(\aleph_0,T)=3$ was found by Peretyatkin in \cite{Peretyat1973complete}. It is obtained by adding an infinite increasing sequence of constants to a dense meet-tree: the axioms state that
$\leqslant$ is a partial order with no maximal elements, for all $x$ the set $\{y\mid y<x\}$ is a dense  linear order without endpoints, that every pair of elements has an $\inf$, and that it is $\omega$-branching: for all $x$ and all $n$ there are distinct $\{y_i\mid i\leqslant n\}$ greater than $x$ such that $\inf(y_i,y_j)=x$ for all $i\neq j$. By adding unary predicates similarly to the case of Ehrenfeucht's example, one can obtain theories with $I(\aleph_0,T)=n$ for all finite $n\geqslant 4$. Here, for any natural number $k\geqslant 2$ we can replace the axioms for $\omega$-branching by the $k$-branching axiom without changing the number $I(\aleph_0,T)$; the $k$-branching axiom states that for all $x$ $k$ is the maximal number such that there are distinct $\{y_i\mid i\leqslant k\}$ greater than $x$ such that $\inf(y_i,y_j)=x$ for all $i\neq j$. Also, in both Ehrenfeucht's and Peretyatkin's examples, we can replace constants by an increasing sequence of unary predicates that define initial parts without sups.

\subsection{Ehrenfeucht theories and SOP} 
Recall that $T$ is said to have the {\em strict order property} (SOP) if there is, in a model of $T$ a definable partial ordering on a definable set with some infinite chains. 
Several examples of Ehrenfeucht theories are known, yet they are essentially the same in the following sense: they were obtained from certain $\omega$-categorical theories with the strict order property, usually dense linear orders or meet trees, by adding axioms for infinitely many constant symbols or unary predicates; each of them has SOP. The following conjecture has circulated since the 1970s and is still widely open. 

\smallskip\noindent
{\bf Conjecture 1.} \  Every Ehrenfeucht theory has the strict order property.  

\smallskip
A weaker (and also open) conjecture is that any (countable, complete) stable theory is not Ehrenfeucht (namely, it has exactly one or infinitely many countable models).  Stability was defined earlier in terms of counting types, but can also be defined as ``not having the order property":  there is no formula $\phi(\bar x,\bar y)$ with $|\bar x|=|\bar y|$ and an infinite sequence of tuples $(\bar a_i\mid i\in\omega)$ from some model of $T$ such that $\models \phi(\bar a_i,\bar a_j)$ if and only if $i<j$. So we see that the strict order property is a special case of stability.

\smallskip 
$T$ is said to have {\em few links} if for any complete types $p(x)$, $q(y)$ of $T$ there are only finitely many complete $r(x,y)$ extending $p(x)\cup q(y)$ which are non-isolated relative to $p$. The notion was introduced by Benda in \cite{Benda}. Pillay in \cite{Pillay1980instability} proved that all Ehrenfeucht theories with few links have the strict order property.
Tsuboi in \cite{Tsuboi1986countable} disproved the expectation that adding an infinite set of constants to some $\aleph_0$-categorical NSOP theory might yield an NSOP Ehrenfeucht theory. He proved that if an Ehrenfeucht theory $T$ is constructed as a countable ascending union of $\omega$-categorical theories, then $T$ has SOP. 
The second author in \cite{Tanovic2006constants} proved that every Ehrenfeucht theory with infinitely many constants has SOP. 

\smallskip   
Lachlan in \cite{lachlan1973number} proved that there is no superstable Ehrenfeucht theory. It was,  on the face of it, a complicated proof. The first author gave a short account of Lachlan's proof, using basic properties of forking,  in \cite{Pillay1983countable}. The key technical point was that superstable theories have a ``rudimentary theory of weight": for $M$ a model of superstable $T$, there do not exist tuples $a$, $b_{1}$, $b_{2}$,..., $b_{i}$,..  such that $\{b_1, b_2, ....\}$ is independent and $a$ forks with each $b_{i}$.  This is a direct consequence of superstability (in terms of its definition that there is no infinite sequence of complete types in variables $x$, $p_{1}(x)\subseteq p_{2}(x \subseteq .....$, each $p_{i+1}$ being a forking extension of $p_{i}$) and basic properties of nonforking independence. 

A rather ``deeper" fact about superstable theories, using regular types, is that all types have ``finite weight", meaning that given $a$ there is a greatest $n$ such that $a$ can fork with each of $b_{1},..,b_{n}$ where $\{b_{1},..,b_{n}\}$ is independent.  And proofs of Lachlan's theorem using finiteness of weight were given by Lascar \cite{Lascar1976ranks} and Shelah \cite{Shelah1990classification}. 

Variants of these properties of weight (and other methods) were used to show that the following classes of theories do not include Ehrenfeuct theories (i.e. $T$ has exactly $1$ or infinitely many countable models).
 
 \begin{itemize}
\item $T$-a countable increasing union of stable theories of finite weight (Tsuboi in \cite{Tsuboi1986countable}).

\item Stable theories admitting finite coding (Hrushovski in \cite{hrushovski1989finitely}).

\item Stable 1-based theories (Pillay \cite{Pillay1989stable}). 

\item Supersimple theories (Kim in \cite{Kim1999number}).

\item NSOP$_1$ theories with non-forking existence and without types of 
self preweight $\omega$  (Kim in \cite{Kim2019number}).

\end{itemize} 

One property of Ehrenfeucht theories, established by Benda in \cite{Benda},  is the existence of powerful types: a non-isolated type $p\in S_n(T)$ is called {\em powerful} if the prime model over its realization realizes all complete types. 
All known examples of powerful types come from known Ehrenfeucht theories; in particular, all of them have SOP and the following is open.

\begin{Problem}\label{P3}
Is there a small NSOP (stable) theory with a powerful type?     
\end{Problem}

In stable Ehrenfeucht theories, powerful types have infinite self preweight. Herwig in \cite{Herwig1995weight} constructed a small stable theory with a type of infinite self preweight; this may suggest the affirmative answer to Problem \ref{P3}.   

\smallskip  
There are a few partial results that confirm Conjecture 1 for theories with 3 countable models.
Woodrow in \cite{Woodrow1978theories} proved that if $T$ is in the same language as the
Ehrenfeucht example, has elimination of quantifiers and 3 countable models,  then T is very much like this example.
Ikeda, Tsuboi and Pillay in \cite{Ikeda1998theories} proved that every almost $\omega$-categorical theory with 3 countable models interprets a dense linear order. ($T$ is almost $\omega$-categorical if for all complete types $p$ and $q$ the type $p(\bar x)\cup q(\bar y)$ has only finitely many completions). 
The second author in \cite{Tanovic2007three} proved that a theory with infinitely many constants and precisely 3 countable models interprets a variant of Ehrenfeucht’s or Peretyatkin’s example.
One may expect that counting the isomorphism types of countable models of binary theories would be relatively easy. 
However, the following problem remains open.

\begin{Problem}
Is there a binary NSOP theory with 3 countable models?   
\end{Problem}

\subsection{Pillay's conjecture} Yet another old open problem on Ehrenfeucht theories was posed by the first author in \cite{Pillay1978number}. 

\smallskip\noindent
{\bf Conjecture 2.} \ If $T$ has a minimal model (i.e. one with no proper elementary substructure), then $T$ has infinitely many countable models.

\smallskip A model $M$ is called algebraic if all its elements are algebraic over $\emptyset$, meaning that if $a\in M$ then there is a formula $\phi(x)$ in the language which is true of $a$
in $M$ and true of only finitely many elements of $M$. 
Note that all algebraic models are minimal, but the converse is not true. A counterexample is $Th(\mathbb Z,S)$, where $S$ is the successor function. However, in \cite{Pillay1978number}, it is shown that if $M$ is prime and minimal, then there exists a formula with parameters $\bar a\in M$ that has infinitely many realizations, all of which are algebraic over $\bar a$. This leads to the following weak version of Conjecture 2, which is still open. 

\smallskip\noindent
{\bf Conjecture 3.} If $T$ has an algebraic model, then $T$ has infinitely many countable models. 

\smallskip
This has been the object of quite a bit of study, which has resulted in proving the following special case.

\smallskip\noindent
{\bf Conjecture 3'.} If $T$ is the complete diagram of a countable structure, then $T$ has infinitely many countable models.

\smallskip
Let $M_0$ be an  algebraic model. In the context of Conjecture 3 (and 3'), by restricting to a formula without parameters of Cantor-Bendixon rank and degree $1$, we may assume that $M_{0}$ itself has CB-rank and degree $1$; we do assume this from now on and let $p(x)$ be the unique nonalgebraic $1$-type of $T$.

In \cite{Pillay1978number} it was proved, in the context of Conjecture 3', that $T$ has at least four countable models.  Note that $T$ cannot be $\omega$-categorical, so by Vaught's result, mentioned earlier, $T$ has at least $3$ countable models, so at least $4$ is the first possible nontrivial statement. This result was extended in \cite{Pillay1980theories} to the case where the elements of $M_0$ are algebraic of uniformly bounded degree (there exists some $k\in\mathbb{N}$ such that for each $a\in M_0$, there exists a $L$-formula satisfied by $a$ that has at most $k$ realizations within $M_0$).

Further progress on Conjecture 3' was made in \cite{Pillay1982dimension}, where it was proved, assuming that the model $M_{0}$ has {\em no order}, namely, we cannot find a formula $\phi(x,y)$  and tuples $a_{1},a_{2},....$  in $M_{0}$ such that $M_{0}\models \phi(a_{i}, a_{j})$ iff $i< j$ for all $i, j\in \omega$.  (This is strictly weaker than saying that $T$ is stable.) Using various bits of machinery developed by the first author, it was shown that one has a well-defined notion of independence among realizations of $p$, and that if $M_{n}$ is prime over $n$-independent realizations of $p$ then the $M_{n}$'s are non isomorphic. Hence $T$ has infinitely many countable models. What is interesting is that much later, the notion of a {\em generically stable type} (first introduced as  ``stable type" by Shelah) was developed in connection with $NIP$ theories (see \cite{Hrushovski-Pillay2011nip}), and then in general in \cite{PillayTanovic2011}, and in fact the type $p$ above was precisely a generically stable, strongly regular type as in \cite{PillayTanovic2011}, and the independence notion was non-forking. 

Finally, the full Conjecture 3' was solved in two papers of the second author, \cite{tanovic2010typesdir} and \cite{Tanovic2011minimal}. In \cite{tanovic2010typesdir}, the notion of a ``type directed by constants over a finite set" is introduced and it is proved that {\em any}  countable complete theory with a type directed by constants has continuum many countable models. 
On the other hand, in \cite{Tanovic2011minimal} a dichotomy theorem was proved for $M_{0}$, $T$ and $p(x)$ as above:  either $p$ is generically stable, strongly regular, or for some infinite $C\subseteq M_{0}$ there is a type directed by $C$ over some finite subset $E$ of the monster model of $T$. Either way we obtain infinitely many countable models. 

\smallskip
Consider again (in the context of Conjecture 3) the algebraic model $M_0$ with CB-rank and degree $1$) and the complete type (over $\emptyset$) $p$ (of CB-rank and degree $1$). Note that the nonalgebraic types in $S_1(M_0)$ correspond precisely to strong types that extend $p$. There are two cases:

\begin{enumerate}[\hspace{10pt} {Case} 1.]
\item $mult (p)=k$ is finite. In particular, this occurs when the elements of $M_0$ are algebraic of uniformly bounded degree. In the subcase when $k=1$, after naming all elements of $M_0$ the dichotomy theorem from \cite{Tanovic2011minimal} applies and Conjecture 3 can be confirmed similarly as above. If $k>1$,
then there is a definable equivalence relation splitting $M_0$ into $k$ infinite classes such that each class with all elements named ''satisfies the assumptions of the dichotomy theorem", so Conjecture 3 follows similarly as before.

\item $mult(p)=2^{\aleph_0}$. In this subcase, Conjecture 3 is wide open. To close it, one may try to extend the dichotomy theorem. Yet another option is trying to prove the following general conjecture. 
\end{enumerate}

 \noindent
{\bf Conjecture 4.} \ Every complete type over $\emptyset$ in an Ehrenfeucht theory has finite multiplicity.

\subsection{Other problems} 
In general, no fine description of the isomorphism types of countable models of Ehrenfeucht theories is known. To support this statement, let us emphasize that the following problem is open.
  
\begin{Problem}\label{Problem_martin_ehrenfeucht}
Prove Martin's conjecture for Ehrenfeucht theories.    
\end{Problem}

This is easy to confirm for theories with 3 countable models by using
Vaught's description of the isomorphism types of their countable models.
Wagner in \cite{CMWagner1982martin} noted that the Scott rank of any Ehrenfeucht theory is strictly below $\omega\cdot\omega$; this bound is still well above $\omega+\omega$ which is predicted by Martin. 

Sudoplatov in \cite{Sudoplatov2004complete} proved that every countable model of a small theory is either ``finitely generated" (prime over a finite set of parameters) or a limit model;  $M$ is a limit model if it can be represented as the union of an elementary $\omega$-chain of finitely generated submodels: $M=\bigcup _{i\in\omega}M_i$, where each $M_i$ is prime over $\bar a_i\in M$. If $T$ is an Ehrenfeucht theory, then the representation can be chosen such that the $M_i$'s are pairwise isomorphic and all tuples $\bar a_i$ realize the same type, say $p\in S_n(T)$; in this case, we say that $M$ is a limit model over $p$. This observation was not enough to resolve  either Problem \ref{Problem_martin_ehrenfeucht} or the following problem, posed by Benda in \cite{Benda}.

\begin{Problem}
Must an Ehrenfeucht theory have a countable, universal, nonsaturated model?   
\end{Problem} 

A (complete) expansion of $T$ by (finitely many) constants is a theory $T_p=T\cup p(\bar c)$, where $p\in S_n(T)$ and $|\bar c|=n$ are new constants; $T_p$'s are also called inessential expansions of $T$.
It is known that by adding constants in this way, we can increase the number of countable models of the theory from 4 to $\aleph_0$ (Woodrow in \cite{Woodrow1978theories}), and even from 3 to $\aleph_0$ (Peretyatkin in \cite{Peretyat1980theories}). There are examples where the number decreases from $\aleph_0$ to finite (Omarov in \cite{Omarov1983nonessential} and Millar in \cite{Millar1989finite}). 
However, it is not possible for the number to increase from $\leqslant \aleph_0$ to $2^{\aleph_0}$: since there are countably many ways to interpret a finite string of new constants in a countable $M\models T$, we have $I(\aleph_0,T_p )\leqslant \aleph_0\cdot I (\aleph_0,T )$. Therefore,
if $T$ has few models, then so does $T_p$. Taimanov in the 1980's asked if the converse is true, but the answer to his question is still not known.   
\begin{Problem}
Is there a theory $T$ with $2^{\aleph_0}$ countable models such that some constant expansion, $T_p$, has $\aleph_0$ (or finitely many) countable models?   
\end{Problem}
Omarov in \cite{Omarov1983nonessential} claimed an affirmative answer to this problem, but, as explained by Baizhanov and Umbetbayev in \cite{BaizhanovUmb2023}, his initial theory had only $\aleph_0$ models. It is interesting that, by adding infinitely many constants (satisfying some fixed type), the number of models can decrease from $2^{\aleph_0}$ to $\aleph_0$ and even to $3$; Simon sketched an example for this in \cite{Hanson}.

\smallskip
As we have already mentioned, Conjecture 1 has been confirmed in \cite{Tanovic2006constants} for theories with infinite $dcl(\emptyset)$; however, for theories with infinite $acl(\emptyset)$, it is still open. In addition, it is natural to consider theories that have a constant expansion with infinite $dcl(\emptyset)$ (or $acl(\emptyset)$).

\begin{Problem}
Must the operators $dcl$ and $acl$ be locally finite in every Ehrenfeucht NSOP theory?     \end{Problem}

This problem is particularly interesting in the context of Ehrenfeucht groups (a group with an additional structure whose theory is Ehrenfeucht). Ehrenfeucht groups exist, an example (group of exponent 2) can be found in \cite{Ilic2015groups}. The following would also be interesting to know. 

\begin{Problem}
 Does every Ehrenfeucht group have  finite exponent?   
\end{Problem}

\printbibliography
\end{document}